\documentclass[12pt]{amsart}
\usepackage{amsfonts}
\usepackage{amsthm}
\usepackage{amsmath}
\usepackage{amscd}
\usepackage[utf8]{inputenc} 
\usepackage{t1enc}
\usepackage[mathscr]{eucal}
\usepackage{indentfirst}
\usepackage{tikz-cd}
\usepackage{hyperref}
\usepackage{graphicx}
\usepackage{graphics}
\usepackage{pict2e}
\usepackage{epic}
\numberwithin{equation}{section}
 \usepackage[T1]{fontenc}
\usepackage{epstopdf} 
\usepackage{tabularx}
\usepackage{xfrac}

\usepackage{xcolor}

\usepackage[tablename=Table]{caption}
\usepackage{relsize}
\usepackage[utf8]{inputenc}
\usepackage{amsmath,amssymb}
\usepackage{amsthm}
\usepackage[all]{xy}
\usepackage{mathrsfs}
\usepackage{tikz-cd}

\newcommand{\ba}{\begin{array}}\newcommand{\ea}{\end{array}}

\usepackage{booktabs,graphicx,makecell,siunitx,eqparbox}

\usepackage{tikz}
\usepackage{tikz-cd}
\usetikzlibrary{%
  matrix,%
  calc,%
  arrows%
}

\usepackage{mathtools}
\usepackage[active]{srcltx}
\usepackage{hyperref,ifpdf}
\usepackage{graphicx}

\newcommand{\Q}{{\mathbb{Q}}}

\newcommand{\C}{{\mathbb{C}}}
\newcommand{\Z}{{\mathbb{Z}}}

\def\C{\mathbb C}

\def\bs#1{\boldsymbol{#1}}

\renewcommand{\dim}{\mathrm{dim}}

\newcommand{\F}{\mathscr F}

\renewcommand{\O}{{\mathcal O}}

\def\P{\mathbb P}
\newcommand{\PP}{\mathbb{P}}

\def\Z{\mathbb Z}

\newcommand{\Sing}{{\rm Sing}}

\newtheorem{lema}{Lemma}[section]

\newtheorem{theorem}{Theorem}
\newtheorem{mcor}{Corollary}

\theoremstyle{definition}
\newtheorem{remark}[lema]{Remark}

\newtheorem{exe}[lema]{Example}

\usepackage{fancyhdr}

\usepackage{lipsum}

\begin{document}

\pagestyle{fancy}


\title{On the connectedness of the singular set of holomorphic foliations}

\author{Omegar Calvo--Andrade}
\address{ Omegar Calvo--Andrade \\ CIMAT: Xalisco SN Valenciana Guanajuato Gto. M\'exico }
\email{jose.calvo@cimat.mx}

\author{Maur\'icio Corr\^ea}
\address{Maur\'icio Corr\^ea  \\ 
Universit\`a degli Studi di Bari, 
Via E. Orabona 4, I-70125, Bari, Italy
}
\email[M. Corr\^ea]{mauricio.barros@uniba.it,mauriciomatufmg@gmail.com } 

\author[  Marcos Jardim]{Marcos Jardim}
\address{IMECC - UNICAMP \\ Departamento de Matem\'atica \\
Rua S\'ergio  Buarque de Holanda, 651\\ 13083-970 Campinas-SP, Brazil}
\email{jardim@ime.unicamp.br}

\author{Jos\'e Seade}
\address{Jos\'e Seade\\ Instituto de Matem\'aticas-Cuernavaca, Universidad Nacional Aut\'onoma de M\'exico.}
\email{jseade@im.unam.mx}

\dedicatory{Dedicated to Professor Nigel Hitchin, with great admiration, \\ on the occasion of his 80th Birthday} 


 \begin{abstract}
Let $\F$ be a singular holomorphic foliation of dimension $k>1$ on a projective $n$-manifold $X$.
Assume that the determinant line bundle $\det(N_{\F})$ of the normal sheaf is ample and that the
singular set $\Sing(\F)$ has dimension $\le k-1$.
We show that the union of the irreducible components of $\Sing(\F)$ of dimension exactly $k-1$
is necessarily connected. This yields a topological obstruction to the
integrability of singular holomorphic distributions, echoing Bott’s vanishing theorem, and it
answers a question by D.~Cerveau for codimension one foliations on $\P^3$.
\end{abstract}

\maketitle

\hyphenation{sin-gu-lari-ties}
\hyphenation{re-si-dues}



 \section*{Introduction}

In the 1960s, Raoul Bott \cite{Bott1970} demonstrated that topological constraints can obstruct the existence of globally integrable foliated structures. In the same vein, we establish a global obstruction to integrability for singular holomorphic distributions. Indeed, by \cite{CCJ, GJM}, there exist holomorphic distributions of dimension $k>1$ on $\P^n$, for all $n \geq 3$, whose singular set contains several connected components of dimension $k-1$. Our main Theorem precludes integrability for all such distributions.

We recall \cite{BB} that a $k$-dimensional \textit{singular} holomorphic foliation   on a complex manifold $X$ of dimension $n$ is defined by  an exact sequence 
$$  0 \to \F \longrightarrow  TX \longrightarrow  N_\F \longrightarrow  0  \,, $$ 
where $TX$ is the sheaf of holomorphic vector fields on $X$, the \textit{tangent sheaf} $\F$ is an integrable coherent subsheaf of $TX$ of 
rank $k$, and the \textit{normal sheaf} $N_\F$ is coherent and  torsion free, of 
rank $n-k;$ this is the codimension of $\F$. The \emph{singular scheme}, denoted by $\Sing(\F)$, is the locus where the normal sheaf is not locally free. This implies that $\Sing(\F)$ has codimension $\ge 2$. Let $\det (N_\F)$ denote the determinant bundle of the normal sheaf. We prove:
 
\begin{theorem}\label{Th. main}
Let $X$ be a projective manifold of dimension $n$,  and $\F$  a singular holomorphic foliation of dimension $k > 1$ such that the bundle $\det (N_\F)$ is ample. Let the singular set be $\Sing(\F) = Z \cup  Z_-$, where $Z$ has dimension $k-1$ and  $ Z_-$ consists of lower-dimensional components. Then $Z$ is a connected set. 
\end{theorem}

Our result is sharp since all the assumptions are necessary; see Section \ref{Examples-remarks}.

Theorem \ref{Th. main} answers positively the question posed by D. Cerveau in \cite{Cerv} regarding dimension 2 foliations on \(\mathbb{P}^3\), since the determinant bundle of the normal sheaf of every holomorphic foliation on projective spaces is ample.

This falls within the Bott phenomenon of global obstructions to integrability: whereas Bott’s theorem compels the vanishing of invariant polynomials in the Chern classes of the normal bundle in degrees exceeding the codimension for integrable subbundles, our connectedness Theorem \ref {Th. main} furnishes a topological obstruction to integrability for singular distributions under an ampleness hypothesis. 

Our proof uses in an essential way the Baum--Bott residues \cite{BB}. The idea is as follows. By Baum--Bott, the Chern number $c_1(N_\F)^{n-k+1}[X]$ has a residue localized at each \textit{connected } component $Z_1, \dots, Z_s$ of the singular set of dimension $k-1$; these residues are elements in $H_{2k-2}(Z_i)$. Then, the ampleness condition yields, after some work,  that the residue at each $Z_i$ equals the total Chern number  $c_1^{n-k+1} (N_\F)\cap [X]$. Hence  $s=1$.

We give examples showing that the ampleness condition on the determinant of the normal sheaf is necessary for Theorem \ref{Th. main}. This condition is always satisfied for holomorphic distributions in $\P^n$, see for instance Remark \ref{remark:Pn}.

 We show in Section \ref{Examples-remarks} that the ampleness hypothesis on $\det(N_{\F})$ is indispensable for Theorem~\ref{Th. main}. 
Moreover, when $X=\PP^n$ and $n\le 2(k-1)$, any two irreducible projective subvarieties $V_1,V_2\subset \PP^n$ of dimension $k-1$ must meet. 
Indeed,  in $A^*(\PP^n)=\mathbb{Z}[h]/(h^{n+1})$  we have  $$
[V_1 \cap V_2]=\deg(V_1)\deg(V_2)h^{\,2n-2k+2}\neq 0. $$
Consequently, the $(k - 1)$-dimensional components of $\Sing(\F)$ cannot be pairwise disjoint, and the conclusion of Theorem~\ref{Th. main} is automatic in this range.
  However, this phenomenon does not hold in general. For example, consider \( X = \mathbb{P}^1 \times \mathbb{P}^5 \) and subvarieties \( W_1 = \{p\} \times \mathbb{P}^4 \) and \( W_2 = \{q\} \times \mathbb{P}^4 \), where \( p \ne q \). Here, \( n = 6 \) and \( k = 4 \), so \( n = 2(k - 1) \), yet \( W_1 \) and \( W_2 \) are clearly disjoint.

In \cite{CCGL}, the first author, Cerveau, Giraldo, and Lins Neto asked whether the tangent sheaf of a codimension one foliation must split as a sum of line bundles whenever it is locally free. Recently, the third author, Faenzi and Vall\'es in \cite{FJMV} provided a negative answer to this problem by producing examples   of rational foliations on \(\mathbb{P}^3\) 
of even degree \(2k+4,~ k\geq 0\), whose tangent sheaf is a twisted null-correlation bundle, hence non-split; the singular sets of such foliations are singular non-reduced curves. Furthermore, these foliations admit flat deformations to non-integrable distributions. 

This naturally leads to the following question: 
\\
\textit{What are the topological obstructions for a rank $2$ reflexive sheaf to be the tangent sheaf of a foliation in 
a projective 3-fold?}

\medskip
We give a partial answer by showing that a natural topological constraint on the tangent sheaf, together with the ampleness of the normal determinant, yields a global obstruction to integrability, detected by the number of isolated singularities of the foliation.

 \begin{mcor}\label{cor:3-folds}
 Let $\F$ be a codimension one foliation on a smooth projective threefold $X$ with $h^1(\mathcal{O}_X)=0$, and let $C$ be the union of the one--dimensional components of $\Sing(\F)$. Assume that $C$ is reduced and that $\det(N_\F)$ is ample. If
\[
h^1(TX\otimes\det(N_\F)^*)=h^2(TX\otimes\det(N_\F)^*)=0,
\]
then
\[
h^2(\F\otimes\det(N_\F)^*) = c_3(\F)\cap [X],
\]
the number of isolated singularities of $\F$ counted with multiplicities. In particular, when $X=\mathbb{P}^3$, these conditions are satisfied, and the formula holds.
\end{mcor}

This is also related to a problem posed by Cerveau in \cite[Problem 2.2]{Cerv}, concerning the description of foliations whose singular set has pure dimension one, equivalently, whose tangent sheaf is locally free (see \cite{CCJ}). Our result provides a topological and algebro-geometric obstruction to the existence of such foliations.

In particular, under the hypotheses of Corollary \ref{cor:3-folds}, a codimension one foliation on a smooth projective threefold with locally free tangent sheaf must satisfy
\[
h^2(\F \otimes \det(N_{\F})^*) = 0,
\]
a vanishing that obstructs a locally free sheaf from being the tangent sheaf of a foliation.

Moreover, a codimension one foliation on $\mathbb{P}^3$ with locally free tangent sheaf forces $\Sing(\F)$ to be singular, unless $\F$ is a pencil of planes.

\begin{mcor}\label{sing:p3}
Let $\F$ be a codimension one foliation on $\mathbb{P}^3$. If the tangent sheaf of $\F$ is locally free, then either $\F$ is a pencil of planes or $\Sing(\F)$ is a singular curve.
\end{mcor}

The hypothesis that $\F$ has locally free tangent sheaf is essential. Indeed, the singular scheme of a generic rational foliation of type $(a,b)$ consists of a smooth complete intersection curve of degree $ab$ together with
\[
(a+b-2)^3+2(a+b-2)^2+2(1-ab)(a+b-2)
\]
distinct points; see \cite[Section 10.1]{CCJ} and \cite{CSV} for the higher-dimensional case. By contrast, in the example of \cite[Theorem 7.1]{FJMV}, where the foliation is given by a pencil and the tangent sheaf is a twisted null-correlation bundle, the singular set is a singular curve.

\subsection*{Acknowledgments}
CA and JS gratefully acknowledge the University of Bari for its hospitality and support through UNIBA Visiting Professor Program.  JS is partially supported by UNAM-PAPIIT grant IN101424. 
MC is partially supported by the Universit\`a degli Studi di Bari and by the PRIN 2022MWPMAB- ``Interactions between Geometric Structures and Function Theories'' and he is a member of INdAM-GNSAGA; 
MJ is supported by the CNPQ grant number 305601/2022-9, the FAPESP-ANR project number 2021/04065-6, and the FAPESP-CEPID project number 2024/00923-6.

We extend our deepest thanks to Nigel J. Hitchin on the occasion of his 80th birthday. His mentorship over the years has been a cornerstone of our mathematical lives, and we are forever grateful for the guidance and inspiration he provided as he showed us the beauty found at the intersection of geometry and topology.

  \section{Proof of Theorem  \ref{Th. main} }
 The proof uses Baum--Bott residues, and by definition, these are elements in the homology groups of $\Sing(\F)$. Yet, it is convenient in the sequel to use these residues  in cohomology, so we consider the following commutative diagram: 
{\small
\[
\begin{tikzcd}[row sep=2em, column sep=2.5em]
& H^{2(n-k+1)}(X, X\setminus Z) \arrow[r, "\phi"] \arrow[dd, "A" pos=0.3] \arrow[dl, "j^*_{U_i}"'] &  H^{2(n-k+1)}(X) \arrow[dd, "P"'] \arrow[dl, "j^*_{U_i}"'] \\
H^{2(n-k+1)}(U_i,  U_i \setminus Z_i) \arrow[r, crossing over, "\phi_i" pos=0.3] \arrow[dd, "A_i"] &  H_c^{2(n-k+1)}(U_i) \arrow[dd,bend right=50,"P_{i}"'] \\
& H_{2k-2}(Z) \arrow[r, "(i_Z)_*"]  & H_{2k-2}(X) \arrow[uu, bend right=50, "P^{-1}"']  \\
H_{2k-2}(Z_i) \arrow[ur, "(i_{Z_i})_*"] 
\arrow[r, "(j_{Z_i})_*"'] 
\arrow[rru,bend  right=35,"(i_{Z_i})_*"'  ] \arrow[uu, bend left=50, "A_{i}^{-1}"] 
&  H_{2k-2}(U_i) \arrow[ur, "(j_{U_i})_*"] \arrow[uu, bend right=50,"P_{i}^{-1}"' pos=0.6]
\end{tikzcd}
\]}
\hskip-4pt
where  all the cohomology and homology groups are with coefficients in 
$\mathbb{C}$; $P=\, \cap [X]$ is    Poincar\'e duality;  the $U_i$ are pairwise disjoint smooth regular neighborhoods (as in \cite{Hirsch}) of the \textit{connected} components $Z_1, \dots,Z_s$  of $Z$, $A$ is Alexander duality,  all the other homomorphisms are induced by the inclusion. The  morphism $\phi$ comes from the following long exact sequence of relative cohomology
\[
\begin{tikzcd}
  H^{2(n-k+1)}(X, X\setminus Z) \arrow[r, "\phi"]  & H^{2(n-k+1)}(X)  \arrow[r] & H^{2(n-k+1)}(X\setminus Z).
\end{tikzcd}
\]
 
Recall that $c_1(N_{\F}) = c_1({\rm det}(N_{\F}))$, see  \cite[Section 5.6]{Kobayashi}.
The residues that we  consider are   localizations of $c_1(N_{\F})^{n-k+1}$ in $H_{2k-2}(Z_i)$. Given one such residue 
\[
{{\rm Res} }_{c_1^{n-k+1}}(\F ,Z_i):={{\rm Res} }(\F ,Z_i)\in H_{2k-2}(Z_i),
\]
we can push it (by the inclusions) to $H_{2k-2}(X)$ and lift it by the Poincar\'e isomorphism to $H^{2(n-k+1)}(X)$. Denote it by  $ {\widehat  {\rm Res} } (\F ,Z_i)$. That is: 
\[
{\widehat  {\rm Res} } (\F ,Z_i) = P^{-1}\big( (i_{(Z_i)})_*  {\rm Res} (\F ,Z_i)\big) \,,
\]
where  ${(i_{(Z_i)_*})} $ denotes the morphism $H_*(Z_i)  \to H_*(Z) \to H_*(X)$. 
Then we have the Baum--Bott formula in the cohomology group $H^{2(n-k+1)}(X)$:
\[
c_1(N_{\F})^{n-k+1} = \sum_{i=1}^s \; {\widehat  {\rm Res} } (\F ,Z_i) \,.
\]

First of all, we observe that the previous formula implies that the ampleness of \( c_1(N_{\F}) \) ensures that the component \( Z = \bigcup_{i=1}^s Z_i \) of dimension \( k-1 \) is non-empty, since \( c_1(N_{\F})^{n-k+1} \neq 0 \), and component $Z_-$ of dimension smaller than \( k-1 \) does not contribute to ${\widehat  {\rm Res} }$.

Now,  consider the above pairwise disjoint smooth regular neighborhoods   of the connected components of $Z$. That is, for each connected component $Z_i\subset Z$, we  take an open neighborhood $U_i$ with a retraction 
\[
\xymatrix{
Z_i  \ar[rr]_{i_{Z_i}} \ar@/^2pc/[rrrr]^{} &&  U_i \ar[rr]_{\rho_i } && Z_i  \\
}
\] 
with a smooth boundary $\partial U_i$, such that 
\[
 \text{for all } j \neq i\quad \overline{U_i} \cap \overline{U_j} = \emptyset, \quad (\overline{U_i}\setminus Z_i)\cap Sing(\F)=\emptyset.
\]

Let $  {\rm Res}  (\F ,Z_i)$ be the Baum--Bott residue in $H_{2k-2}(Z_i)$. Define:
\[
 \widehat {\rm Res}  (\F  |_{U_i},Z_i) = P_i^{-1}\big((i_{Z_i})_* {\rm Res} (\F ,Z_i)\big)  \,.
\]

By the commutativity of the diagram above, we have the following equality in $H_c^{2(n-k+1)}(U_i)$:
\begin{equation}\label {eq 1}
  \widehat {\rm Res}  (\F  |_{U_i},Z_i) = j^*_{U_i} \big( \widehat  {\rm Res}  (\F ,Z_i) \big) \;.
\end{equation}

We can   take a neighborhood $V_i$ of  $ \partial U_i $ such that $\F$ is  regular  on $V_i$, i.e.,  the tangent sheaf of $\F$ is locally free on $V_i$, where     $U_i$   is  a retraction of $Z_i$. Then, we can apply Theorem 1.4 from \cite{SS1} and, using the naturality of Chern classes, we obtain
\begin{equation} \label{eq 2}
\begin{aligned}
\widehat{\mathrm{Res}}(\F|_{U_i}, Z_i) &= c_1(\det(N_{\F})|_{U_i})^{n-k+1} \\
\\
&= c_1(j^*_{U_i} \det(N_{\F}))^{n-k+1} \\
\\
&= j^*_{U_i}\big(c_1(N_{\F})^{n-k+1}\big) \;.
\end{aligned}
\end{equation}

Write $Z_i$ as the union of its irreducible components of dimension $k-1$:
\[
Z_i=\bigcup_{\alpha=1}^{m_i} W_{i\alpha}\qquad(\dim W_{i\alpha}=k-1).
\]
The localized residue in $U_i$ then reads
\[
\widehat{\mathrm{Res}(\F|_{U_i},Z_i)}
=\sum_{\alpha=1}^{m_i}\lambda_{i\alpha}\,[\widehat{W_{i\alpha}}]
\ \in H_c^{2(n-k+1)}(U_i).
\]
Applying $j^*_{U_i}$ we get the following equality in $H_c^{2(n-k+1)}(U_i)$:
\begin{equation}\label{eq 4}
j^*_{U_i} \big(\widehat {\rm Res}(\F ,Z_i)\big) \,=\, \sum_{\alpha=1}^{m_i}\lambda_{i\alpha}\, j^*_{U_i}\big([\widehat{W_{i\alpha}}]\big)  \,=\,  j^*_{U_i}  \big( c_1({\rm det}(N_{\F})\big)^{n-k+1} \big),
\end{equation}
where the last equality follows from equation (\ref{eq 2}).

Set $L:=\det(N_{\F})$ and write $a:=c_1(L)\in H^2(X;\Q)$.
Choose $m\ge 1$ such that $L^{\otimes m}$ is very ample and set
$$h:=c_1(L^{\otimes m})\in H^2(X;\Z)\subset H^2(X;\Q).$$
Then $a=(1/m )\cdot h$ in $H^2(X;\Q)$.
Since $k>1$, cup \eqref{eq 4} with $h^{k-1}$, viewed in $H^{2(k-1)}(X;\Q)$ and then pulled back to $H_c^{2(k-1)}(U_i;\Q)$).
We obtain the following equality  in $H_c^{2n}(U_i;\Q)$: 
\[
\sum_{\alpha=1}^{m_i}\lambda_{i\alpha}\, j^*_{U_i}\!\big([\widehat{W_{i\alpha}}]\smile h^{k-1}\big)
= j^*_{U_i}\!\big(c_1(L)^{n-k+1}\smile h^{k-1}\big).
\]
For each $\alpha$, define
\[
d_{i\alpha}:=\int_X [W_{i\alpha}]\smile h^{k-1}\in\Z_{>0}.
\]
Equivalently, if $H_1^i,\dots,H_{k-1}^i\in |L^{\otimes m}|$ are sufficiently general and
\[
W_i:=H_1^i\cap\cdots\cap H_{k-1}^i,
\]
then $W_{i\alpha}\cap W_i$ is a $0$--dimensional scheme and
$
d_{i\alpha}=\mathrm{length}(W_{i\alpha}\cap W_i).
$
Thus
\[
[\widehat{W_{i\alpha}}]\smile h^{k-1}=d_{i\alpha}\,\eta
\qquad\text{in }H^{2n}(X;\Q),
\]
where $\eta$ is a fixed generator of the one-dimensional space $H^{2n}(X;\Q)$.
Observe that we do not identify $\eta$ with $h^n$; we only use $\eta\neq 0$.
Therefore, after evaluating in top degree, we get
\[
\sum_{\alpha=1}^{m_i}\lambda_{i\alpha}\, d_{i\alpha}
=\int_X c_1(L)^{n-k+1}\smile h^{k-1}
=\frac{1}{m^{n-k+1}}\int_X h^{n}\neq 0.
\]
In particular, the same nonzero number appears for every $i$.
On the other hand, Baum--Bott gives
\[
c_1(L)^{n-k+1}=\sum_{i=1}^{s}\widehat{\mathrm{Res}}(\F,Z_i)
=\sum_{i=1}^{s}\sum_{\alpha=1}^{m_i}\lambda_{i\alpha}\,[\widehat{W_{i\alpha}}]
\ \text{in }H^{2(n-k+1)}(X;\Q).
\]
Cupping with $h^{k-1}$ and integrating over $X$ yields
\[
\int_X c_1(L)^{n-k+1}\smile h^{k-1}
=\sum_{i=1}^{s}\sum_{\alpha=1}^{m_i}\lambda_{i\alpha}\, d_{i\alpha}
=s\cdot \int_X c_1(L)^{n-k+1}\smile h^{k-1}.
\]
Since $$\int_X c_1(L)^{n-k+1}\smile h^{k-1}\neq 0,$$ we conclude $s=1$.
Hence $Z$ is connected.

\qed

 \section{Proof of Corollary 1}
 
  Theorem \ref{Th. main}  guarantees that the curve $C$   is connected. Applying \cite[Theorem 5]{CCJ2}, we conclude that $$h^2(\F\otimes\det(N_\F)^*)=c_3(\F)\cap [X],$$ and by 
\cite[Proposition 3.6]{Gholampour-Kool}
$$
c_3(\F)\cap [X]= h^0(X,\mathscr{E}xt^1(\F,\mathcal{O}_X) )
$$
which is the number of isolated singularities of $\F$ counted with multiplicities, since by \cite[Lemma 2.1]{CCJ}, the foliation $\F$ has tangent sheaf locally free along $C$. 
Reducedness alone does not exclude embedded associated points in general.
In what follows we only use the elementary fact that for a reduced curve $C$,
$h^0(\mathcal O_C)$ equals the number of connected components of $C$.
In particular, if $C$ is reduced and connected, then $h^0(\mathcal O_C)=1$.

Now, suppose \( X = \mathbb{P}^3 \). Then, the dual of the determinant of the normal sheaf \( N_{\mathscr{F}} \) is given by 
$
\mathcal{O}_{\mathbb{P}^3}(-2 - d),
$
where \( d \) denotes the degree of the foliation \( \mathscr{F} \), see Remark \ref{remark:Pn}. Applying Bott's formulae \cite[pg.12]{AMM}, we deduce
\[
h^1(\mathcal{O}_{\mathbb{P}^3})  = h^1(T\mathbb{P}^3(-2 - d))\simeq h^1(\Omega_{\mathbb{P}^3}^1(2 - d))  = 0 \quad \text{for all } d \geq 0.
\]
Furthermore, following the proof of Proposition 14 in \cite{CCJ2}, we obtain
\[
h^2(T\mathbb{P}^3 \otimes \det(N_{\mathscr{F}})^*) = h^2(T\mathbb{P}^3(-2 - d)) = h^1(\mathscr{I}_C) = h^0(\mathcal{O}_C) - 1,
\]
which vanishes whenever \( C \) is a reduced and connected curve.

\qed

  \section{Proof of Corollary 2}
According to \cite[Theorem 3.11]{CCJ}, a codimension one distribution $\F$ of positive degree on $\mathbb{P}^3$ with a locally free tangent sheaf and a smooth, connected singular scheme cannot be integrable. Remark \ref{remark:Pn} below states that Theorem \ref{Th. main} applies to this situation. Therefore, if \( \F \) is integrable, then either \( \F \) has degree zero and is a pencil of planes in which case \( \Sing(\F) \) is a line, or \( \Sing(\F) \) cannot be smooth.

\section{Examples and remarks}\label{Examples-remarks}

In this section we give examples illustrating the necessity of the hypotheses and the sharpness of Theorem~1.

\begin{exe}
Let $X$ be the classical Hopf manifold of dimension $n$. This is obtained from $\C^n \setminus \{0\}$ dividing it by the group of automorphisms generated by a contraction, say $z \mapsto z/2 $. It is clear that  $X$ fibers over $\P^{n-1}$ with fiber an elliptic curve, the quotient of $\C^*$ by the contraction. So one has 
an elliptic fibration $X  \buildrel {\pi} \over {\to} \P^{n-1}$.

Consider now a 1-dimensional foliation $\F$  on $\P^{n-1}$ with isolated singularities. Its pull back  $\pi^{\ast}\F$ is a 2-dimensional foliation  on $X$ with 
 singular set, a union of disjoint elliptic curves:
\[
\operatorname{Sing}(\pi^*\F) = \bigsqcup_{p \in \operatorname{Sing}(\F)} \pi^{-1}(p).
\]
 If $\operatorname{Sing}(\F)$ is reduced, then $\operatorname{Sing}(\pi^*\F)$ has $d^{n-1}+ \cdots + d+1$ connected components, see for instance \cite[Theorem A]{Soares}.

In this example, neither the manifold in question is projective, nor is $\det(N_\F)$ ample. 
\end{exe}

\begin{exe}(Foliations on  products) 
Let $\F$ be a holomorphic foliation of dimension $1$ with isolated singularities on a projective manifold $X$  of dimension $n$. If $Y$ is a projective manifold of dimension $k-1$, then we have a foliation $\pi_2^*\F$ on $X\times Y$ of dimension $k$,  with singular set the  disjoint union 
\[
\operatorname{Sing}(\pi_2^*\F) = \bigsqcup_{p \in \operatorname{Sing}(\F)} \pi_2^{-1}(p),
\]
where $\pi_2^{-1}(p)=Y$ for all $p \in \operatorname{Sing}(\F)$. Then $$c_1(\det(N(\pi_2^*\F)))=\pi_2^*c_1(\det(N_\F)) ,$$ and 
\[
c_1(\det(N_\F(\pi_2^*\F)))^{n+k-1}=0
\]
since $n+k-1>n$.  So the determinant of the normal bundle is not ample. 
\end{exe}

The next example is a 2-dimensional foliation on \(\mathbb{P}^n\), induced by the action of the affine group, whose singular set is the connected union of \(n\) rational curves intersecting at a single point.

\begin{exe} This example generalizes to higher dimensions the example of \cite{CCGL} for $n=3$. The foliation $\widetilde{\F}_{3,d}$ below defines a rigid irreducible component of codimension one foliations on $\P^3$ of degree $d+1$, and $\widetilde{\F}_{3,1}$ generates the \emph{exceptional component} of degree $2$ foliations \cite[Cor.~3, p.~1007]{CCGL}.

In the affine space $\{[1:z]\}\simeq\C^n\subset\P^n,\, n\geq3$, let us denote $c_k(d)=c_k(T\P^k(d-1))$, and we define the complete vector fields 
\[
\xymatrix@R=1mm{ u_{n,d}(z) = z_1\partial_{z_1}+c_1(d)\cdot z_2\partial_{z_2}+\cdots+ c_{n-1}(d)\cdot z_n\partial_{z_n} \\
v_{n,d}(z)= \partial_{z_1}+c_1(d)\cdot z_1^d\partial_{z_2}+\cdots+ c_{n-1}(d)\cdot 
z_{n-1}^d\partial_{z_n}. } 
\]  

Since  
$[v_{n,d},u_{n,d}]=v_{n,d},$ they are the infinitesimal generators of an action of the affine group, $ \mathfrak{aff}(\mathbb{C})\times \C^n\longrightarrow\C^n.$

The vector field $u_{n,d}$ is linear, it extends as a holomorphic section $\widetilde{u}_{n,d}$ of $T\P^n$. On the other hand,  
$v_{n,d}$ is polynomial of degree $d$, it extends as a holomorphic bundle map $\widetilde{v}_{n,d}:\O_{\P^n}(1-d)\longrightarrow T\P^n$, with a unique singularity at the point $p=[0:\cdots:0:1]\in \P^n$, with Milnor number
\[
\mu(\widetilde{v}_{n,d}(p))=c_n(T\P^n(d-1))=c_{n}(d).
\]
We define the foliation of dimension 2:
\[
\xymatrix{ 0\ar[r] & \widetilde{\F}_{n,d}=\O_{\P^n}\oplus\O_{\P^n}(1-d) 
\ar[rr]^-{\widetilde{u}_{n,d}\oplus\widetilde{v}_{n,d} } && T\P^n \ar[r] & N_{\widetilde{\F}_{n,d}}\ar[r] &0   }
\]
The singular set $\Sing(\widetilde{\F}_{n,d})=\Lambda_{1,d}\cup\Lambda_{2,d}\cup\dots\cup\Lambda_{n,d}$ is a union of rational curves
$\Lambda_{j,d}$ of degrees $\deg(\Lambda_{1,d})=1$ and for $j>1,\quad \deg(\Lambda_{j,d})=c_{j-1}(d).$ Moreover 
\[
\bigcap_{j=1}^n \Lambda_{j,d}=\{ [0:\cdots:0:1] \},
\]
hence, $\Sing(\widetilde{\F}_{n,d})$ is connected. 
Notice that the restriction $\widetilde{\F}_{n,d}|_{\{Z_n=0\}}$ is $\widetilde{\F}_{n-1,d}$, and 
\[
 \Sing(\widetilde{\F}_{n,d})=\Sing(\widetilde{\F}_{n-1,d})\cup \Lambda_{n,d}.
\] 
Denote $C_k=c_k(d)$, so  the curve $\Lambda_{n,d}$ is parameterized by 
\[
\Lambda_{n,d}([s:t])=[ s^{C_{n-1}}:s^{(C_{n-1}-1)}t:s^{(C_{n-1}-C_1)}t^{C_1}:\cdots :t^{C_{n-1}}].
\]  
\end{exe}

\begin{remark}
Let $\F$ be a holomorphic foliation on $\PP^n$ of dimension $k>1$, and let
$Z\subset \Sing(\F)$ denote the (non-empty) union of the irreducible components
of $\Sing(\F)$ of pure dimension $k-1$. Write
\[
Z=\bigcup_{j=1}^{s} Z_j,
\]
where the $Z_j$ are the connected components of $Z$.
Assume that $n\le 2k-2$. Then any two projective subvarieties of $\PP^n$ of
dimension $k-1$ must meet. Indeed, for each $j$ we have
$[Z_j]=d_j\,h^{n-k+1}\in H^{2(n-k+1)}(\PP^n)$, where $h\in H^2(\PP^n)$ is the
hyperplane class and $d_j=\deg(Z_j)>0$. Hence, for any indices $j\neq \ell$,
\[
\begin{aligned}
[Z_j]\smile[Z_\ell]
&= (d_j h^{n-k+1})\smile(d_\ell h^{n-k+1}) \\
&= d_j d_\ell\,h^{2(n-k+1)}\in H^{4(n-k+1)}(\PP^n),
\end{aligned}
\]
and this class is nonzero because $2(n-k+1)\le n$. Then, $Z_j\cap Z_\ell\neq\emptyset$
for all $j\neq \ell$, and in particular $Z$ is connected. Consequently, the content of
Theorem  \ref{Th. main} is only nontrivial in the range $n>2k-2$.
Finally, we note that in the range $n\ge 2k$, it was proved in \cite{CorreaSuwa2025}
that $\dim ( \Sing(\F))\ge k-1$, so $Z$ is nonempty in that case.
\end{remark}

\begin{remark}\label{exe-non-conn}
We note that the hypothesis that $k>1$ is also necessary: a generic foliation by curves on a complex projective manifold has 0-dimensional singularities consisting of several distinct points. In addition, \cite{CJM} contains various examples of foliations by curves on $\mathbb{P}^3$ whose singular set is the union of two or more skew lines.
\end{remark}

\begin{remark}\label{remark:Pn}
Consider a holomorphic foliation \(\F \subset TX\) of dimension \(k\).  
By taking the maximal exterior power of the dual morphism 
\[
\Omega^1_X \longrightarrow \F^*,
\] 
one obtains a morphism
$
\wedge^k \phi^\vee: \Omega^k_X \longrightarrow \det(\F)^*.
$
This, in turn, yields an induced morphism
$$
\Omega^k_X \otimes \det(\F) \longrightarrow \mathcal{O}_X,
$$
whose image is an ideal sheaf \(\mathcal{I}_Z\) defining the  singular set of the foliation \(Z=\Sing(\F) \subset X\).  
Moreover, this construction determines a global section
$
\nu_{\F} \in H^0\left(X, \wedge^k TX \otimes \det(\F)^*\right). 
$
By the isomorphism 
$
\wedge^k TX \simeq \Omega_X^{n-k} \otimes \det(TX)
$
and the relation 
$
\det(N_\F) \simeq \det(TX) \otimes K\F,
$
we conclude that the isomorphism
\[
\wedge^k TX \otimes K\F 
\simeq \Omega_X^{n-k} \otimes \det(TX) \otimes K\F 
\simeq \Omega_X^{n-k} \otimes \det(N_\F)
\]
shows that the section \(\nu_{\F}\) corresponds to a twisted \((n-k)\)-form
\[
\omega_{\F} \in H^0\left(X, \Omega_X^{n-k} \otimes \det(N_\F)\right)
\]
with values in \(\det(N_\F)\).
A holomorphic foliation \(\F\) of dimension \(k\) on \(\mathbb{P}^n\) can be described by a twisted \((n-k)\)-form
\[
\omega \in H^0\left(\mathbb{P}^n, \Omega_{\mathbb{P}^n}^{n-k}(d+n-k+1)\right),
\]
where \(d \geq 0\) is referred to as the \emph{degree} of the foliation \(\F\). Then, \(\det(N_\F)=\mathcal{O}_{\mathbb{P}^n}(d+n-k+1)\)   is   ample. See for instance  \cite[Section 2.5]{CJV}.
\end{remark}

\begin{remark}
Observe that, a priori, each connected component is such that 
\[
Z_i = \bigcup_{j=1}^{j(i)} Z_{ij}
\]
where \( Z_{ij} \) are its irreducible components, for $i=1,\dots,s$.  
Theorem \ref{Th. main} implies  the component of dimension \( k-1 \) of the singular set of the foliation \( Z \) is connected. That is,    $Z=Z_i$, for all $i$,  and  by Baum--Bott
$$c_1(N_{\F})^{n-k+1} \cap [X] = \sum_{j=1}^\ell \lambda _j  [Z_j] \in H_{2k-2}(Z),$$
where \( Z_{j} \) are  the irreducible components of $Z$ and \(\lambda_j\) are complex numbers that depend only on the local behavior of the leaves of \(\F\) near \(Z_j\).
\end{remark}

\end{document}